\documentclass[11pt]{article}
\usepackage{mathrsfs}
\usepackage{epic,eepic,epsf,epsfig}
\usepackage{amsfonts,srcltx,mathrsfs}
\usepackage{amsmath}
\usepackage{amssymb}
\usepackage{amsbsy}
\usepackage{graphicx}
\usepackage{amsfonts}
\usepackage{color}
\usepackage{mathrsfs}
\usepackage{epic,eepic,epsf,epsfig}
\usepackage{graphicx}
\usepackage{amsfonts,srcltx,mathrsfs}
\usepackage{array}
\usepackage{amsthm}
\usepackage{algorithm}
\usepackage{algorithmic}
\newtheorem{theorem}{Theorem}
\newtheorem{lem}{Lemma}

\newtheorem{defi}{Definition}

\newtheorem{clm}{Claim}

\def\f{\noindent}

\begin{document}

\markboth{ et. al}{The generalized connectivity of $(n,k)$-bubble-sort graphs}

\title{The generalized connectivity of $(n,k)$-bubble-sort graphs
}

\author{Shu-Li Zhao$^{1}$, Rong-Xia Hao$^{1}$\footnote{Corresponding author. Email:17118434@bjtu.edu.cn(S.-L. Zhao), rxhao@bjtu.edu.cn (R.-X. Hao), lw@uttyler.edu(L. Wu)},  Lidong Wu$^{2}$\\[0.2cm]
{\em\small $^{1}$Department of Mathematics, Beijing Jiaotong University,}\\ {\small\em
Beijing 100044, P.R. China}\\ {\em\small $^{2}$Department of Computer Science, University of Texas at Tyler,}\\ {\small\em Tyler, Texas 75799, USA}}

\date{}
\maketitle

Let $S\subseteq V(G)$ and $\kappa_{G}(S)$ denote the maximum number $r$ of edge-disjoint trees $T_{1}, T_{2}, \cdots,\\ T_{r}$ in $G$ such that $V(T_{i})\bigcap V(T_{j})=S$ for any $i, j \in \{1, 2, \cdots, r\}$ and $i\neq j$. For an integer $k$ with $2\leq k\leq n$, the {\em generalized $k$-connectivity} of a graph $G$ is defined as $\kappa_{k}(G)= min\{\kappa_{G}(S)|S\subseteq V(G)$ and $|S|=k\}$. The generalized $k$-connectivity is a generalization of the traditional connectivity. In this paper, the generalized $3$-connectivity of the $(n,k)$-bubble-sort graph $B_{n,k}$ is studied for $2\leq k\leq n-1$. By proposing an algorithm to construct $n-1$ internally disjoint paths in $B_{n-1,k-1}$, we show that $\kappa_{3}(B_{n,k})=n-2$ for $2\leq k\leq n-1$, which generalizes the known result about the bubble-sort graph $B_{n}$ [Applied Mathematics and Computation 274 (2016) 41-46] given by Li $et$ $al.$, as the bubble-sort graph $B_{n}$ is the special $(n,k)$-bubble-sort graph for $k=n-1$.

\medskip

\f {\em Keywords:} Generalized connectivity; fault-tolerance; interconnection network; $(n,k)$-bubble-sort graph.

\section{Introduction}
For an interconnection network, one mainly concerns about the reliability and fault tolerance. An interconnection network is usually modelled as a connected graph $G=(V, E)$, where nodes represent processors
and edges represent communication links between processors.  {\it The connectivity $\kappa (G)$} of a graph $G$ is an important parameter to evaluate the reliability and fault tolerance of a network. It is defined as the minimum number of vertices whose deletion results in a disconnected graph. In addition, Whitney~\cite{w} defines the connectivity from local point of view. That is, for any subset $S=\{u, v\}\subseteq V(G)$, let $\kappa_{G}(S)$ denote the maximum number of internally disjoint paths between $u$ and $v$ in $G$. Then $\kappa(G)=min\{\kappa_{G}(S)|S\subseteq V(G)$ and $|S|=2\}.$ As a generalization of the traditional connectivity, Chartrand et al.~\cite{c} introduced the {\it generalized $k$-connectivity} in $1984$. This parameter can measure the reliability of a network $G$ to connect any $k$ vertices in $G$. Let $S\subseteq V(G)$ and $\kappa_{G}(S)$ denote the maximum number $r$ of edge-disjoint trees $T_{1}, T_{2}, \ldots, T_{r}$ in $G$ such that $V(T_{i})\bigcap V(T_{j})=S$ for any $i, j \in \{1, 2, \ldots, r\}$ and $i\neq j$. For an integer $k$ with $2\leq k\leq n$, the {\it generalized $k$-connectivity} of a graph $G$ is defined as $\kappa_{k}(G)= min\{\kappa_{G}(S)|S\subseteq V(G)$ and $|S|=k\}$. The generalized $2$-connectivity is exactly the traditional connectivity. Li~\cite{l4} derived that it is NP-complete for a general graph $G$ to decide whether there are $k$ internally disjoint trees connecting $S$, where $k$ is a fixed integer and $S\subseteq V(G).$ Some results~\cite{l2,l5} about the upper and lower bounds of the generalized connectivity are obtained. In addition, there are some results of the generalized $k$-connectivity for some classes of graphs and most of them are about $k=3$. For example, Chartrand {\em et al.}~\cite{ch} studied the generalized connectivity of complete graphs; Li {\em et al.}~\cite{LIS} characterized the minimally $2$-connected graphs with generalized connectivity $\kappa_{3}=2$; Li {\em et al.}~\cite{l1} studied the generalized $3$-connectivity of Cartesian product graphs; Li {\em et al.}~\cite{l8} studied the generalized $3$-connectivity of graph products; Li {\em et al.}~\cite{l3} studied the generalized connectivity of the complete bipartite graphs; Li {\em et al.}~\cite{l6} studied the generalized $3$-connectivity of the star graphs and bubble-sort graphs; Li {\em et al.}~\cite{l7} studied the generalized $3$-connectivity of the Cayley graph generated by trees and cycles and Lin and Zhang~\cite{L} studied the generalized $4$-connectivity of hypercubes etc.

In this paper, we focus on the $(n,k)$-bubble-sort graph, denoted by $B_{n,k}$. The complete graph $K_{n}$ and the bubble-sort graph $B_{n}$ are special $(n,k)$-bubble-sort graphs $B_{n,k}$ for $k=1$ and $k=n-1$, respectively. In~\cite{ch}, it was shown that $\kappa_{3}(K_{n})=n-2$ for $n\geq 3$ and in~\cite{l6}, it was shown that $\kappa_{3}(B_{n})=n-2$ for $n\geq 3$. Following, we study the generalized $3$-connectivity of $B_{n,k}$ for $2\leq k\leq n-1$ and it is shown that $\kappa_{3}(B_{n,k})=n-2$, which generalizes the known results about bubble-sort graphs~\cite{l6}.

The paper is organized as follows. In section 2, some notation and definitions are given. In section 3, the connectivity of $(n,k)$-bubble-sort graphs $B_{n,k}$ is determined for $2\leq k\leq n-1$. In addition, the generalized $3$-connectivity of $B_{n,k}$ is determined for $2\leq k\leq n-1$ and an algorithm for constructing $n-1$ internally disjoint paths in $B_{n-1,k-1}$ was proposed. In section 4, the paper is concluded.

\section{Preliminary}

Let $G=(V, E)$ be a simple, undirected graph. Let $|V(G)|$ be the size of vertex set and $|E(G)|$ be the size of edge set. For a vertex $v$ in $G$, we denote by $N_{G}(v)$ the {\em neighbourhood} of the vertex $v$ in $G$ and $N_{G}[v]=N_{G}(v)\bigcup\{v\}$. Let $U \subseteq V(G)$, denote $N_G(U)=\bigcup\limits_{v\in U}N_{G}(v)-U$. Let $d_{G}(v)$ denote the degree of the vertex $v$ in $G$ and $\delta(G)$ denote the {\em minimum degree} of the graph $G$. The subgraph induced by $V^{\prime}$ in $G$, denoted by $G[V^{\prime}]$, is a graph whose vertex set is $V^{\prime}$ and the edge set is the set of all the edges of $G$ with both ends in $V^{\prime}$. A graph is said to be {\em $k$-regular} if for any vertex $v$ of $G$, $d_{G}(v)=k$. Two $xy$- paths $P$ and $Q$ in $G$ are {\em internally disjoint} if they have no common internal vertices, that is $V(P)\bigcap V(Q)=\{x, y\}$. Let $Y\subseteq V(G)$ and $X\subset V(G)\setminus Y$, the $(X, Y)$-paths is a family of internally disjoint paths starting at a vertex $x\in X$, ending at a vertex $y\in Y$ and whose internal vertices belong neither to $X$ nor $Y$. If $X=\{x\}$, the $(X, Y)$-paths is a family of internal disjoint paths whose starting vertex is $x$ and the terminal vertices are distinct in $Y$, which is referred to as a {\em $k$-fan} from $x$ to $Y$. For terminologies and notation not undefined here we follow the reference~\cite{B}.

Let $\Gamma$ be a finite group and $S$ be a subset of $\Gamma$, where the identity of the group does not belong to $S$. The {\em Cayley graph $Cay(\Gamma, S)$} is a digraph with vertex set $\Gamma$ and arc set $\{(g, g.s)| g\in \Gamma, s\in S\}$. If $S= S^{-1}$, then $Cay(\Gamma, S)$ is an undirected graph, where $S^{-1}=\{s^{-1}|s \in S\}$.

Let $[n]=\{1,2,\cdots,n\}$ and $Sym(n)$ denote the group of all permutations on $[n]$. Let $(p_{1}p_{2}\cdots p_{n})$ denote a permutation on $[n]$ and $(ij)$, which is called a transposition, denote the transposition that swaps the objects at positions $i$ and $j$, that is, $(p_{1}\cdots p_{i}\cdots p_{j}\cdots p_{n})(ij)=(p_{1}\cdots p_{j}\cdots p_{i}\cdots p_{n})$. For the Cayley graph $Cay(Sym(n), T)$, where $T$ is a set of transpositions of $Sym(n)$. Let $G(T)$ be the graph on $n$ vertices $\{1,2,\ldots,n\}$ such that there is an edge $ij$ in $G(T)$ if and only if transposition $(ij)\in T$~\cite{s}. The graph $G(T)$ is called {\em the transposition generating graph} of $Cay(Sym(n), T)$. It is well known that if $G(T)\cong P_{n}$, where $P_{n}$ is a path with $n$ vertices, then $Cay(Sym(n), T)$ is called an {\em $n$-dimensional bubble sort graph} and denoted by $B_{n}$.

As a generalization of $B_{n}$, the $(n,k)$-bubble-sort graph, denoted by $B_{n,k}$, was introduced by Shawash~\cite{Sha} in $2008$. The $(n,k)$-bubble-sort graph $B_{n,k}$ is defined as follows.

\begin{defi}\label{defi2}
Given two positive integers $n$ and $k$ with $n>k$, let $[n]$ denote the set $\{1, 2, \cdots, n\}$ and $P_{n,k}$ be a set of arrangements of $k$ elements in $[n]$. The $(n,k)$-bubble-sort graph $B_{n,k}$ has vertex set $P_{n,k}$, and two vertices $u=a_{1}a_{2}\cdots a_{k}$ and $v=b_{1}b_{2}\cdots b_{k}$ are adjacent if and only if one of the following conditions hold.
\begin{enumerate}
\item [{\rm (a)}] There exists an integer $m\in [2, k]$ such that $a_{m-1}=b_{m}, a_{m}=b_{m-1}$ and $a_{i}=b_{i}$ for all $i\in[k]\setminus \{m-1, m\}$.

\item [{\rm (b)}] $a_{i}=b_{i}$ for all $i\in[k]\setminus \{1\}$ and $a_{1}\neq b_{1}$.
\end{enumerate}
\end{defi}

For two distinct $i$ and $j$, where $i\in [n]$ and $j\in [k]$. Let $V_{n,k}^{j:i}$ be the set of vertices in $B_{n,k}$ with the $j$th position being $i$, that is, $V_{n,k}^{j:i}=\{p|p=p_{1}p_{2}\cdots p_{j}\cdots p_{k}\in P_{n,k}$ and $p_{j}=i\}$. For a vertex $v=p_{1}p_{2}\cdots p_{i}\cdots p_{n}$, we call $p_{i}$ the element at position $i$ of the vertex $v$. For a fixed position $j\in[k]$, $\{V_{n,k}^{j:i}|1\leq i \leq n\}$ forms a partition of $V_{n,k}$. Let $B_{n,k}^{j:i}$ denote the subgraph of $B_{n,k}$ induced by $V_{n,k}^{j:i}$. Then for each $j\in [k]$, $B_{n,k}^{j:i}$ is isomorphic to $B_{n-1, k-1}$. Thus, $B_{n,k}$ can be recursively constructed from $n$ copies of $B_{n-1, k-1}$. It is easy to check that each $B_{n,k}^{j:i}$ is a subgraph of $B_{n,k}$ and $B_{n,k}$ can be decomposed into $n$ subgraphs $B_{n,k}^{j:i}$s according to the $j$th position. By the symmetry of $B_{n,k}$ and for simplicity, we shall take $j$ as the last position $k$ and use $B_{n,k}^{i}$ to denote $B_{n,k}^{k:i}$. For convenience, let $B_{n,k}=B_{n,k}^{1}\bigoplus B_{n,k}^{2}\bigoplus \cdots \bigoplus B_{n,k}^{n}$, where $\bigoplus$ just denotes the corresponding decomposition of $B_{n,k}$. Obviously, any vertex $u$ of $B_{n,k}^{i}$ has $k-1$ neighbors in $B_{n,k}^{i}$ and one neighbor outside of $B_{n,k}^{i}$, which is called the outside neighbour of $u$.
\begin{figure}[!ht]
\begin{center}
\vskip1cm
\includegraphics[scale=0.2]{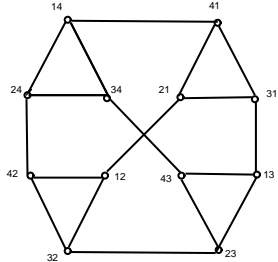}
\end{center}
\vskip0.5cm
\caption{The $(4, 2)$-bubble-sort graph $B_{4,2}$}\label{F2}
\end{figure}
Let $E(i,j)$ be the set of edges between $B_{n,k}^{i}$ and $B_{n,k}^{j}$, that is, $E(i,j)=\{(p, q)\in E(B_{n,k})|p\in V(B_{n,k}^{i})$ and $q\in V(B_{n,k}^{j})\}$. Clearly, $E(i,j)$ is a matching between $B_{n,k}^{i}$ and $B_{n,k}^{j}$ and $|E(i,j)|=\frac{(n-2)!}{(n-k)!}$. By the definition of $B_{n,k}$, $B_{n,1}$ is isomorphic to $K_{n}$ and $B_{n,n-1}$ is isomorphic to $B_{n}$. It follows that $B_{n,k}$ is a generalization of the bubble-sort graph $B_{n}$. The $(4, 2)$-bubble-sort graph $B_{4,2}$ is depicted in Figure~\ref{F2}.

\section{The generalized $3$-connectivity  of the $(n, k)$-bubble-sort graph }

In this section, the generalized $3$-connectivity of the $(n, k)$-bubble-sort graph $B_{n,k}$ will be proved. To prove the result, the following lemmas are useful.

\begin{lem}\label{lem1}
Let $B_{n, k}=B_{n, k}^{1}\bigoplus B_{n, k}^{2}\bigoplus \ldots\bigoplus B_{n, k}^{n}$ for $n\geq k+1$ and $1\leq k\leq n-1$. Then the following results hold.
\begin{enumerate}
\item [{\rm (1)}] For any vertex $u$ of $B_{n, k}^{i}$, it has exactly one outside neighbour.

\item [{\rm (2)}] For any copy $B_{n, k}^{i}$, no two vertices in $B_{n, k}^{i}$ have a common outside neighbour. In addition, $|N(B_{n,k}^{i})|=\frac{(n-1)!}{(n-k)!}$ and $|N(B_{n, k}^{i})\bigcap V(B_{n, k}^{j})|=\frac{(n-2)!}{(n-k)!}$ for $i\neq j$.
\end{enumerate}
\end{lem}

\f {\bf Proof.} (1) By the definition of $B_{n, k}$, the result holds clearly.

(2) Let $u, v\in V(B_{n,k}^{i})$ and $u\neq v$. If they have a common outside neighbour $w$, then $u$ and $v$ are the two outside neighbours of $w$ which lie in the same copy, which contradicts with {\rm (1)}. Thus, no two vertices in $B_{n, k}^{i}$ have a common outside neighbour.

Since $|V(B_{n,k}^{i})|=\frac{(n-1)!}{(n-k)!}$ and no two vertices in $B_{n,k}^{i}$ have a common outside neighbor, $|N(B_{n,k}^{i})|=\frac{(n-1)!}{(n-k)!}$ and $|N(B_{n, k}^{i})\bigcap V(B_{n, k}^{j})|=\frac{(n-2)!}{(n-k)!}$ for $i\neq j$.
\hfill\qed

\begin{lem}{\rm(\cite{l5})}\label{lem2}
Let $G$ be a connected graph and $\delta$ be its minimum degree. Then $\kappa_{3}(G)\leq \delta$. Further, if there are two adjacent vertices of degree $\delta$, then $\kappa_{3}(G)\leq \delta-1$.
\end{lem}

\begin{lem}{\rm(\cite{l5})}\label{lem3}
Let $G$ be a connected graph with $n$ vertices. If $\kappa(G)=4k+r,$ where $k$ and $r$ are two integers with $k\geq 0$ and $r\in \{0, 1, 2, 3\},$ then $\kappa_{3}(G)\geq 3k+\lceil\frac{r}{2}\rceil$. Moreover, the lower bound is sharp.
\end{lem}

\begin{lem}{\rm(\cite{B})}\label{lem4}
Let $G=(V, E)$ be a $k$-connected graph, and let $X$ and $Y$ be subsets of $V(G)$ of cardinality at least $k$. Then there exists a family of $k$ pairwise disjoint $(X, Y)$-paths in $G$.
\end{lem}

\begin{lem}{\rm(\cite{B})}\label{lem5}
Let $G=(V, E)$ be a $k$-connected graph, let $x$ be a vertex of $G$, and let $Y\subseteq V\setminus \{x\}$ be a set of at least $k$ vertices of $G$. Then there exists a $k$-fan in $G$ from $x$ to $Y$, that is, there exists a family of $k$ internally disjoint $(x, Y)$-paths whose terminal vertices are distinct in $Y$.
\end{lem}

Next, we determine the connectivity of $B_{n,k}$ for $k=2$.

\begin{lem}\label{lem6}
$\kappa(B_{n,2})=n-1$ for $n\geq 3$.
\end{lem}

\f {\bf Proof.} Let $B_{n, 2}=B_{n, 2}^{1}\bigoplus B_{n, 2}^{2}\bigoplus \ldots\bigoplus B_{n, 2}^{n}$. Let $F$ be a minimum vertex cut of $B_{n, 2}$ and $u\in V(B_{n, 2})$. Since $N_{B_{n, 2}}(u)$ is a vertex cut of $B_{n, 2}$ and $|N_{B_{n, 2}}(u)|=n-1$, $|F|\leq n-1$.

Next, we show that $|F|\geq n-1$. Suppose to the contrary, that is, $|F|\leq n-2$. Let $F_{i}=F\bigcap V(B_{n, 2}^{i})$ for each $i\in\{1, 2,\cdots,n\}$. Without loss of generality, let $|F_{1}|\geq |F_{2}|\geq\cdots\geq |F_{n}|$. Then $|F_{_{n-1}}|=|F_{_{n}}|=0$. By Lemma~\ref{lem1}(2), $B_{n, 2}[V(B_{n, 2}^{n-1})\bigcup V(B_{n, 2}^{n})]$ is connected. Let $C$ be a component of $B_{n, 2}-F$ that does not contain $B_{n, 2}[V(B_{n, 2}^{n-1})\bigcup V(B_{n, 2}^{n})]$ as a subgraph and $c_{i}=|V(C)\bigcap V(B_{n, 2}^{i})|$ for each $i\in\{1, 2,\cdots, n-2\}$. Then there exists an integer $l\in\{1, 2,\cdots, n-2\}$ such that $c_{l}> 0$. Let $u\in V(B_{n, 2}^{l})\bigcap V(C)$ and $u^{\prime}\in V(B_{n, 2}^{j})$, where $u^{\prime}$ is the outside neighbour of $u$ in $B_{n, 2}^{j}$, $j\in [n]$ and $l\neq j$.

If $u^{\prime}\in V(B_{n, 2}^{j})\setminus V(C)$, then $u^{\prime}\in F_{j}$. It implies that $|F_{j}|\geq 1$.

If $u^{\prime}\in V(C)$, then $N_{B_{n, 2}^{j}}(V(B_{n, 2}^{n-1})\bigcup V(B_{n, 2}^{n}))\subseteq F_{j}$. Otherwise, the component that contains $B_{n, 2}[V(B_{n, 2}^{n-1})\bigcup V(B_{n, 2}^{n})]$ will be $C$ as $B_{n, 2}^{j}\cong K_{n-1}$, which is a contradiction. By Lemma~\ref{lem2}, $|N_{B_{n, 2}^{j}}(V(B_{n, 2}^{n-1})\bigcup V(B_{n, 2}^{n})|=2$. It implies that $|F_{j}|\geq 2$.

Recall that $B_{n, 2}^{l}$ is a complete graph, then $|F|=|F_{1}\bigcup \cdots \bigcup F_{n}|\geq |V(B_{n, 2}^{l})|-c_{l}+c_{l}=n-1$, a contradiction. Thus, $|F|\geq n-1$.

\hfill\qed

Next, we determine the connectivity of $B_{n, k}$ for $2\leq k\leq n-1$.

\begin{lem}\label{lem7}
$\kappa(B_{n,k})=n-1$ for $2\leq k\leq n-1$.
\end{lem}
\f {\bf Proof.} Let $F$ be a minimum vertex cut of $B_{n, k}$ and $u\in V(B_{n, 2})$. Since $N_{B_{n, k}}(u)$ is a vertex cut of $B_{n, k}$ and $|N_{B_{n, k}}(u)|=n-1$, $|F|\leq n-1$.

Next, we show that $\kappa(B_{n,k})\geq n-1$. We prove the result by induction on $k$. When $n\geq 3$ and $k=2$, by Lemma~\ref{lem6}, the result holds. Suppose that the result holds for $B_{n^{\prime},k-1}$, where $2\leq k-1\leq n^{\prime}-2$. Now we consider $B_{n,k}$ for $3\leq k\leq n-2$. Let $F_{i}=F\bigcap V(B_{n, k}^{i})$ for each $i\in\{1, 2,\cdots,n\}$. Without loss of generality, let $|F_{1}|\geq |F_{2}|\geq\cdots\geq |F_{n}|$. Suppose to the contrary, that is, $|F|\leq n-2$. Thus, $|F_{n-1}|=|F_{n}|=0$.

If $|F_{1}|=n-2$, then $|F_{i}|=0$ for each $i\in\{2,3,\cdots,n\}$. By Lemma~\ref{lem1}(2), $B_{n, k}[\bigcup_{i=2}^{n}V(B_{n, k}^{i})]$ is connected. As any vertex in $B_{n, k}^{1}\setminus F_{1}$ has an outside neighbour, $B_{n, k}-F$ is connected, a contradiction.

If $|F_{1}|\leq n-3$, then $|F_{i}|\leq n-3$ for each $i\in\{2,3,\cdots,n\}$. By induction, $B_{n, k}^{i}-F_{i}$ is connected for each $i\in\{1,2,\cdots,n\}$. As $|F_{n}|=0$ and there are $\frac{(n-2)!}{(n-k)!}$ independent edges between $B_{n, k}^{i}$ and $B_{n, k}^{n}$. Note that $\frac{(n-2)!}{(n-k)!}-|F_{i}|\geq \frac{(n-2)!}{(n-3)!}-|F_{i}|\geq 1$ for each $i\in\{1,2,\cdots,n-1\}$. Then there exists at least one edge between $B_{n, k}^{i}-F_{i}$ and $B_{n, k}^{n}$. It implies that $B_{n, k}-F$ is connected, a contradiction. Thus, $|F|\geq n-1$.

\hfill\qed

To prove the main result, the following lemmas are useful.

\begin{lem}\label{lem8}
Let $B_{n, k}=B_{n, k}^{1}\bigoplus B_{n, k}^{2}\bigoplus \ldots\bigoplus B_{n, k}^{n}$ and $H=B_{n, k}[V(B_{n, k})\setminus V(B_{n, k}^{i})]$ for some $i\in[n]$. If $2\leq k\leq n-1$, then $\kappa(H)=n-2$.
\end{lem}

\f {\bf Proof.} Without loss of generality, let $H=B_{n, k}[V(B_{n, k})\setminus V(B_{n, k}^{n})]$, that is, $H=B_{n, k}^{1}\bigoplus B_{n, k}^{2}\\\bigoplus \ldots\bigoplus B_{n, k}^{n-1}$. As there is some vertex $v\in V(H)$ whose outside neighbour belongs to $B_{n, k}^{n}$, $\delta(H)=n-2$. Hence, $\kappa(H)\leq \delta (H)=n-2$.

Next, we show that $\kappa(H)\geq n-2$. To prove the result, we just need to show that for any two distinct vertices $v_{1}$ and $v_{2}$ of $H$, there exist at least $n-2$ internally disjoint paths between them. The result is proved by considering the following two cases.

Case 1. $v_{1}$ and $v_{2}$ belong to the same copy of $B_{n-1,k-1}$.

Without loss of generality, let $v_{1}, v_{2}\in V(B_{n, k}^{1})$. By Lemma~\ref{lem7}, $\kappa(B_{n, k}^{1})=n-2$. Hence, there are $n-2$ internally disjoint paths between $v_{1}$ and $v_{2}$ in $B_{n, k}^{1}$.

Case 2. $v_{1}$ and $v_{2}$ belong to different copies of $B_{n-1,k-1}$.

Without loss of generality, let $v_{1}\in V(B_{n, k}^{1})$ and $v_{2}\in V(B_{n, k}^{2})$.

Subcase 2.1. $3\leq k\leq n-1$

By Lemma~\ref{lem1}(2), there are $\frac{(n-2)!}{(n-k)!}$ independent edges between $B_{n, k}^{1}$ and $B_{n, k}^{2}$. Choose $n-2$ vertices $u_{1}, u_{2}, u_{3},\cdots,$ $u_{n-2}$ from $B_{n, k}^{1}$ such that the outside neighbour $u_{i}^{\prime}$ of $u_{i}$ belongs to $B_{n, k}^{2}$ for each $i\in\{1,2,\cdots,n-2\}$. This can be done as $\frac{(n-2)!}{(n-k)!}\geq n-2$ for $k\geq 3$ and $n\geq k+1$. Let $S=\{u_{1}, u_{2}, u_{3},\cdots, u_{n-2}\}$ and $S^{\prime}=\{u_{1}^{\prime}, u_{2}^{\prime}, u_{3}^{\prime},\cdots, u_{n-2}^{\prime}\}$. By Lemma~\ref{lem7}, $\kappa(B_{n, k}^{1})=\kappa(B_{n, k}^{2})=n-2$. If $v_{1}\notin S$, by Lemma~\ref{lem5}, there exists a family of $n-2$ internally disjoint $(v_{1}, S)$-paths $P_{1}, P_{2}, \cdots, P_{n-2}$ whose terminal vertices are distinct in $S$. Note that if $v_{1}\in S$, then there is a $(v_{1}, S)$ path that contains the only vertex $v_{1}$. Similarly, if $v_{2}\notin S^{\prime}$, by Lemma~\ref{lem5}, there exists a family of $n-2$ internally disjoint $(v_{2}, S^{\prime})$ paths $P_{1}^{\prime}, P_{2}^{\prime}, \cdots, P_{n-2}^{\prime}$ whose terminal vertices are distinct in $S^{\prime}$. Note that if $v_{2}\in S^{\prime}$, there is a $(v_{2}, S^{\prime})$ path that contains the only vertex $v_{2}$. Let $\widehat{P_{i}}=P_{i}\bigcup u_{i}u_{i}^{\prime}\bigcup P_{i}^{\prime}$ for each $i\in\{1,2,\cdots,n-2\}$, then $n-2$ disjoint paths between $v_{1}$ and $v_{2}$ are obtained in $H$.

Subcase 2.2. $k=2$ and $n\geq 3$

By Lemma~\ref{lem1}(2), there is exactly one edge between $B_{n, k}^{i}$ and $B_{n, k}^{j}$ for $i\neq j$ and $i,j\in\{1,2,\cdots,n-1\}$. Choose $n-2$ vertices $u_{1}, u_{2}, u_{3},\cdots,$ $u_{n-2}$ from $B_{n, k}^{1}$ such that the outside neighbour $u_{i}^{\prime}$ of $u_{i}$ belongs to $B_{n, k}^{i+1}$ for each $i\in\{1,2,\cdots,n-2\}$, and choose $n-3$ vertices $w_{2}, w_{3},\cdots,$ $w_{n-2}$ from $B_{n, k}^{2}$ such that the outside neighbour $w_{i}^{\prime}$ of $w_{i}$ belongs to $B_{n, k}^{i+1}$ for each $i\in\{2,3,\cdots,n-2\}$. Let $S=\{u_{1}, u_{2}, u_{3},\cdots, u_{n-2}\}$ and $S^{\prime}=\{u_{1}^{\prime}, w_{2}, w_{3},\cdots, w_{n-2}\}$. Note that $B_{n, k}^{i}\cong K_{n-1}$ for each $i\in\{1,2,\cdots,n\}$. If $v_{1}\notin S$, then $S=N_{B_{n, k}^{1}}(v_{1})$. If $v_{1}\in S$, let $v_{1}=u_{1}$. Then $S\setminus \{u_{1}\}\subseteq N_{B_{n, k}^{1}}(v_{1})$. Similarly, if $v_{2}\notin S^{\prime}$, then $S^{\prime}=N_{B_{n, k}^{2}}(v_{2})$. If $v_{2}\in S^{\prime}$, let $v_{2}=u_{1}^{\prime}$. Then $S^{\prime}\setminus \{u_{1}^{\prime}\}\subseteq N_{B_{n, k}^{2}}(v_{2})$. Recall that $B_{n, k}^{i}\cong K_{n-1}$ for $i\in[n-1]$, then $u_{i}^{\prime}w_{i}^{\prime}$ is an edge in $B_{n, k}^{i+1}$ for each $i\in\{2,3,\cdots,n-2\}$. Let $P_{1}=v_{1}u_{1}u_{1}^{\prime}v_{2}$ and $P_{i}=v_{1}u_{i}u_{i}^{\prime}w_{i}^{\prime}w_{i}v_{2}$ for each $2\leq i\leq n-2$, then $n-2$ disjoint paths between $v_{1}$ and $v_{2}$ are obtained in $H$.

Hence, $\kappa(H)=n-2$.
\hfill\qed

\begin{lem}\label{lem9}
Let $B_{n, 2}=B_{n, 2}^{1}\bigoplus B_{n, 2}^{2}\bigoplus \ldots\bigoplus B_{n, 2}^{n}$. For any vertex $v\in V(B_{n, 2}^{i})$ for $1\leq i\leq n$, let $N_{B_{n, 2}^{i}}[v]=N_{B_{n, 2}^{i}}(v)\bigcup\{v\}$. Then $|N_{B_{n, 2}^{i}}[v]|=n-1$ and the $n-1$ outside neighbours of vertices in $N_{B_{n, 2}^{i}}[v]$ belong to different copies of $B_{n-1, 1}$.
\end{lem}

\f {\bf Proof.} Let $v\in V(B_{n, 2}^{i})$, then $d_{B_{n, 2}^{i}}(v)=n-2$. Thus, $|N_{B_{n, 2}^{i}}[v]|=n-1$ holds clearly. Without loss of generality, assume $i=2$ and $v=12$. Then $N_{B_{n, 2}^{i}}[v]=\{32, 42, \cdots,n2\}$. Let $S$ be the set of  outside neighbours of the vertices in $N_{B_{n, 2}^{i}}[v]$, then $S=\{21, 23, 24, \cdots, 2n\}$. Hence, the outside neighbours are contained in $B_{n, 2}^{1}, B_{n, 2}^{3}, \cdots, B_{n, 2}^{n}$, respectively. The result is desired.
\hfill\qed

Following, we prove the generalized $3$-connectivity of $B_{n,k}$ for $k=2$.

\begin{theorem}\label{thm1}
$\kappa_{3}(B_{n,2})=n-2$ for $n\geq 3$.
\end{theorem}

\f {\bf Proof.} As $B_{n,2}$ is $(n-1)$-regular. By Lemma~\ref{lem2}, $\kappa_{3}(B_{n,2})\leq \delta-1= n-2$. To complete the result, it suffices to show that $\kappa_{3}(B_{n,2})\geq n-2$. We prove the result by induction on $n$.

For $n=3$, $B_{3,2}$ is connected. Then $\kappa_{3}(B_{3,2})\geq 1=n-2$.

For $n=4$, by Lemma~\ref{lem3} and Lemma~\ref{lem7}, $\kappa_{3}(B_{n, 2})\geq \lceil\frac{3}{2} \rceil=2=n-2$.

Next, suppose that $n\geq 5$. Let $B_{n, 2}=B_{n, 2}^{1}\bigoplus B_{n, 2}^{2}$ $\bigoplus\ldots \bigoplus B_{n,2}^{n}$ and $v_{1}, v_{2}, v_{3}$ be any three distinct vertices of $B_{n,2}$. For convenience, let $S=\{v_{1}, v_{2}, v_{3}\}$. We prove the result by considering the following three cases.

Case 1. $v_{1}, v_{2}$ and $v_{3}$ belong to the same copy of $B_{n-1,1}$.

Without loss of generality, let $v_{1}, v_{2}, v_{3}\in V(B_{n,2}^{1})$. By the inductive hypothesis, $\kappa_{3}(B_{n,2}^{1})\\\geq n-3$. That is, there are $n-3$ internally disjoint trees $T_{1}, T_{2}\cdots, T_{n-3}$ connecting $S$ in $B_{n,2}^{1}$. Let $v_{1}^{\prime}, v_{2}^{\prime}$ and $v_{3}^{\prime}$ be the outside neighbours of $v_{1}, v_{2}$ and $v_{3}$, respectively. Then $\{v_{1}^{\prime}, v_{2}^{\prime}, v_{3}^{\prime}\}\subseteq V(B_{n, 2})\setminus V(B_{n, 2}^{1})$. As $B_{n, 2}[V(B_{n, 2})\setminus V(B_{n, 2}^{1})]$ is connected, there exists a tree $T$ connecting $v_{1}^{\prime}, v_{2}^{\prime}$ and $v_{3}^{\prime}$ in $B_{n, 2}[V(B_{n, 2})\setminus V(B_{n, 2}^{1})]$. Let $T_{n-2}=T\bigcup v_{1}v_{1}^{\prime}\bigcup v_{2}v_{2}^{\prime}\bigcup v_{3}v_{3}^{\prime}$, then it is a tree connecting $S$ and $V(T_{n-2})\bigcap V(B_{n,2}^{1})=S$. Hence, there exist $n-2$ internally disjoint trees connecting $S$ in $B_{n,2}$ and the result is desired.

Case 2. $v_{1}, v_{2}$ and $v_{3}$ belong to two different copies of $B_{n-1, 1}$.

Without loss of generality, let $v_{1}, v_{2}\in V(B_{n,2}^{1})$ and $v_{3}\in V(B_{n,2}^{2})$. By Lemma~\ref{lem7}, $\kappa(B_{n,2}^{1})=n-2$. Hence, there exist $n-2$ internally disjoint paths $P_{1}, P_{2}, \ldots, P_{n-2}$ between $v_{1}$ and $v_{2}$ in $B_{n,2}^{1}$. Choose $n-2$ distinct vertices $x_{1}, x_{2}, \ldots, x_{n-2}$ from $P_{1}, P_{2}, \ldots, P_{n-2}$ such that $x_{i}\in V(P_{i})$ for each $i\in \{1,2,\cdots,n-2\}$. Note that at most one of these paths has length $1$. If there is one path with length $1$, say $P_{1}$ and let $x_{1}=v_{1}$. Let $x_{i}^{\prime}$ be the outside neighbour of $x_{i}$ for each $i\in \{1,2,\cdots,n-2\}$. Let $X^{\prime}=\{x_{1}^{\prime}, x_{2}^{\prime}, \cdots, x_{n-2}^{\prime}\}$, then $X^{\prime}\subset V(B_{n,2})\setminus V(B_{n,2}^{1})$. By Lemma~\ref{lem1}, $|X^{\prime}|=n-2$. By Lemma~\ref{lem8}, $B_{n,2}[V(B_{n,2})\setminus V(B_{n,2}^{1})]$ is $n-2$ connected. By Lemma~\ref{lem5}, there exist $n-2$ internally disjoint $(v_{3}, X^{\prime})$-paths $P_{1}^{\prime}, P_{2}^{\prime}, \ldots, P_{n-2}^{\prime}$ in $B_{n,2}[V(B_{n,2})\setminus V(B_{n,2}^{1})]$ whose terminal vertices are distinct in $X^{\prime}$. Note that if $v_{3}\in X^{\prime}$, then there is a $(v_{3}, X^{\prime})$-path that contains exactly one vertex $v_{3}$. Let $T_{i}=P_{i}\bigcup x_{i}x_{i}^{\prime}\bigcup P_{i}^{\prime}$ for each $i\in\{1,2,\cdots,n-2\}$. Then $n-2$ internally disjoint trees connecting $S$ in $B_{n,2}$ are obtained.

Case 3. $v_{1}, v_{2}$ and $v_{3}$ belong to three different copies of $B_{n-1,1}$, respectively.

Without loss of generality, let $v_{1} \in V(B_{n,2}^{1}), v_{2} \in V(B_{n,2}^{2})$ and $v_{3} \in V(B_{n,2}^{3})$. Let $N_{B_{n,2}^{i}}[v_{i}]=N_{B_{n,2}^{i}}(v_{i})\bigcup\{v_{i}\}$ for $i=1,2,3$. By Lemma~\ref{lem9}, for each $i\in \{1,2,3\}$ and $j\in \{4,5,\cdots, n\}$, there exists one vertex in $N_{B_{n, 2}^{i}}[v_{i}]$, say $u_{i}^{j}$, such that the outside neighbour $(u_{i}^{j})^{\prime}$ of $u_{i}^{j}$ belongs to $B_{n,2}^{j}$. As $B_{n,2}^{j}$ is connected, we can find a tree $\widehat{T}_{j}$ connecting $(u_{1}^{j})^{\prime}, (u_{2}^{j})^{\prime}$ and $(u_{3}^{j})^{\prime}$ for each $j\in\{4,5,\cdots,n\}$.
Let $T_{j}=\widehat{T}_{j}\bigcup u_{1}^{j}(u_{1}^{j})^{\prime} \bigcup u_{2}^{j}(u_{2}^{j})^{\prime} \bigcup u_{3}^{j}(u_{3}^{j})^{\prime}\bigcup v_{1}u_{1}^{j}\bigcup v_{2}u_{2}^{j}\\\bigcup v_{3}u_{3}^{j}$ as $B_{n-1,1}\cong K_{n-1}$, then $n-3$ internally disjoint trees connecting $S$ are obtained. Let $\widehat{B}_{n,2}^{i}=B_{n,2}^{i}-(\{u_{i}^{4}, u_{i}^{5},\cdots, u_{i}^{n}\}\setminus \{v_{i}\})$. Then there are at most $n-3$ vertices deleted from $B_{n,2}^{i}$ for each $i\in\{1,2,3\}$. As $B_{n,2}^{i}$ is $n-2$ connected, $\widehat{B}_{n,2}^{i}$ is still connected. For $i,j \in \{1,2,3\}$ and $i\neq j$, there is exactly an edge between $B_{n,2}^{i}$ and $B_{n,2}^{j}$.  Thus, $B_{n,2}[\bigcup_{i=1}^{3}V(\widehat{B}_{n,2}^{i})]$ is connected and there is a tree $T_{n-2}$ connecting $S$. Hence, there exist $n-2$ internally disjoint trees connecting $S$ in $B_{n,2}$ and the result is desired.
\hfill\qed

Next, we prove the generalized $3$-connectivity of $B_{n,k}$ for $3\leq k\leq n-1$.

\begin{theorem}\label{thm2}
$\kappa_{3}(B_{n,k})=n-2$ for $3\leq k\leq n-1$.
\end{theorem}

\f {\bf Proof.} As $B_{n,k}$ is $(n-1)$-regular. By Lemma~\ref{lem2}, $\kappa_{3}(B_{n,k})\leq \delta-1= n-2$. To complete the result, it suffices to show that $\kappa_{3}(B_{n,k})\geq n-2$. We prove the result by induction on $n$.

For $n=3$, $B_{3,k}$ is connected. Then $\kappa_{3}(B_{3,k})\geq 1=n-2$.

For $n=4$, by Lemma~\ref{lem3} and Lemma~\ref{lem7}, $\kappa_{3}(B_{n,k})\geq \lceil\frac{3}{2} \rceil=2=n-2$.

Next, suppose that $n\geq 5$. Let $B_{n,k}=B_{n,k}^{1}\bigoplus B_{n,k}^{2}$ $\bigoplus\ldots \bigoplus B_{n,k}^{n}$ and $v_{1}, v_{2}, v_{3}$ be any three distinct vertices of $B_{n,k}$. For convenience, let $S=\{v_{1}, v_{2}, v_{3}\}$. We prove the result by considering the following three cases.

Case 1. $v_{1}, v_{2}$ and $v_{3}$ belong to the same copy of $B_{n-1,k-1}$.


Case 2. $v_{1}, v_{2}$ and $v_{3}$ belong to two different copies of $B_{n-1,k-1}$.


Case 3. $v_{1}, v_{2}$ and $v_{3}$ belong to three different copies of $B_{n-1,k-1}$, respectively.

The proofs of Case $1$ and Case $2$ are the same as the proof of Case $1$ and Case $2$ in Theorem $1$. Thus, only the Case $3$ is considered.

Without loss of generality, let $v_{1} \in V(B_{n,k}^{1}), v_{2} \in V(B_{n,k}^{2})$ and $v_{3} \in V(B_{n,k}^{3})$. Let $v_{1}=p_{1}p_{2}\cdots p_{k-1}1$ and $v_{i}=p_{i}p_{2}\cdots p_{k-1}1$ for $k+1\leq i\leq n$, where $p_{k+1}, p_{k+2}, \cdots, p_{n}$ are distinct elements in $[n]\setminus \{p_{1}, p_{2},\cdots, p_{k-1},1\}$. We now present the algorithm, called (n-1)IDP, that constructs $n-1$ internally disjoint paths $P_{2}^{1}, P_{3}^{1},\cdots, P_{n}^{1} $ in $B_{n}^{1}$ such that the outside neighbour of each terminal vertex of the $n-1$ paths belong to different copies of $B_{n-1, k-1}$.

\begin{algorithm}[h]
\caption{(n-1)IDP(k) }
\begin{algorithmic}[1]
\REQUIRE $n, k$, where $3\leq k \leq n-1$, $v_{1}=p_{1}p_{2}\cdots p_{k-1}1$;
\ENSURE $n-1$ pairwise disjoint path $P_{2}^{1}, P_{3}^{1},\cdots,P_{k}^{1}, P_{k+1}^{1},\cdots, P_{n}^{1}$;
\FOR{$i=2$ to $k-1$}
\STATE $P_{i}^{1}=v_{1}, t=v_{1}$;\
\FOR{$j=i$ to $k-1$}
\STATE $t=t(j-1,j)$ // where $(j-1,j)$ is a transposition
\STATE $P_{i}^{1}=P_{i}^{1}\bigcup t;$\
\ENDFOR
\ENDFOR
\STATE $P_{k}^{1}=v_{1}$;\
\FOR{$i=k+1$ to $n$}
\STATE $P_{i}^{1}=v_{1}v_{i}$, $t=v_{i}=p_{i}p_{2}\cdots p_{k-1}1$;\
\FOR{$j=1$ to $k-2$}
\STATE $t=t(j,j+1)$ // where $(j,j+1)$ is a transposition
\STATE $P_{i}^{1}=P_{i}^{1}\bigcup t;$\
\ENDFOR
\ENDFOR
\end{algorithmic}
\end{algorithm}

By the above algorithm, there are the following $n-1$ paths $P_{2}^{1}, P_{3}^{1},\cdots, P_{n}^{1}$ starting at the vertex $v_{1}$ in $B_{n,k}^{1}$, where $p_{k+1}, p_{k+2}, \cdots, p_{n}$ are distinct elements in $[n]\setminus \{p_{1}, p_{2},\cdots, p_{k-1},1\}$.

$P_{2}^{1}=(\underline{p_{1}}p_{2}p_{3}\cdots p_{k-1}1)(i_{2}\underline{p_{1}}p_{3}\cdots p_{k-1}1)(p_{2}p_{3}\underline{p_{1}}\cdots p_{k-1}1)\cdots (p_{2}p_{3}\cdots p_{k-1}\underline{p_{1}}1)$;

$P_{3}^{1}=(p_{1}\underline{p_{2}}p_{3}\cdots p_{k-1}1)(p_{1}p_{3}\underline{p_{2}}\cdots p_{k-1}1)\cdots (p_{1}p_{3}\cdots p_{k-1}\underline{p_{2}}1)$;

$\cdots $

$P_{k-1}^{1}=(p_{1}p_{2}p_{3}\cdots \underline{p_{k-2}}p_{k-1}1)(p_{1}p_{2}p_{3}\cdots p_{k-1}\underline{p_{k-2}}1)$;

$P_{k}^{1}=(p_{1}p_{2}p_{3}\cdots \underline{p_{k-1}}1)$;

$P_{k+1}^{1}=(p_{1}p_{2}p_{3}\cdots p_{k-1}1)(\underline{p_{k+1}}p_{2}p_{3}\cdots p_{k-1}1)(p_{2}\underline{p_{k+1}}p_{3}\cdots p_{k-1}1)(p_{2}p_{3}\underline{p_{k+1}}\cdots p_{k-1}1)\cdots(p_{2}\\p_{3}p_{4}\cdots\underline{p_{k+1}}1)$;

$P_{k+2}^{1}=(p_{1}p_{2}p_{3}\cdots p_{k-1}1)(\underline{p_{k+2}}p_{2}p_{3}\cdots p_{k-1}1)(p_{2}\underline{p_{k+2}}p_{3}\cdots p_{k-1}1)(p_{2}p_{3}\underline{p_{k+2}}\cdots p_{k-1}1)\cdots(p_{2}\\p_{3}p_{4}\cdots \underline{p_{k+2}}1)$;

$\cdots $

$P_{n}^{1}=(p_{1}p_{2}p_{3}\cdots p_{k-1}1)(\underline{p_{n}}p_{2}p_{3}\cdots p_{k-1}1)(p_{2}\underline{p_{n}}p_{3}\cdots p_{k-1}1)(p_{2}p_{3}\underline{p_{n}}\cdots p_{k-1}1)\cdots(p_{2}p_{3}p_{4}\cdots\underline{p_{n}}\\1)$.

\begin{clm}\label{clm1}
For every $a,b\in\{2,3,\cdots,n\}$ and $a\neq b, V(P_{a}^{1})\bigcap V(P_{b}^{1})=\{v_{1}\}$.
\end{clm}

The proof of the Claim $1$. Without loss of generality, suppose that $a< b$.

If $a,b\in\{2,3,\cdots,k\}$, then for any vertex $y\in V(P_{a}^{1})\setminus\{v_{1}\}$, the $a-1$ elements at positions $1,2,\cdots,a-1$ of $y$ are $p_{1}p_{2}\cdots p_{a-2}p_{a}$. However, for any vertex $z\in V(P_{b}^{1})\setminus\{v_{1}\}$, the $a-1$ elements at positions $1,2,\cdots,a-1$ of $z$ are $p_{1}p_{2}\cdots p_{a-2}p_{a-1}$. As $p_{a}\neq p_{a-1}$, then $y\neq z$. Hence, the claim holds.

If $a, b\in\{k+1,\cdots,n\}$, then for any vertex $y\in V(P_{a}^{1})\setminus\{v_{1}\}$, it is the permutation of $\{p_{a},p_{2}\cdots,p_{k-1},1\}$. For any vertex $z\in V(P_{b}^{1})\setminus\{v_{1}\}$, it is the permutation of $\{p_{b},p_{2}\cdots,p_{k-1},1\}$. As $p_{a}, p_{b}\in [n]\setminus \{p_{1},p_{2}\cdots,p_{k-1},1\}$ and $p_{a}\neq p_{b}$, then $y\neq z$. Thus, the claim holds.

If $a\in\{2,3,\cdots,k\}$ and $b\in\{k+1,\cdots,n\}$, then for any vertex $y\in V(P_{a}^{1})\setminus\{v_{1}\}$, it is the permutation of $\{p_{1},p_{2}\cdots,p_{k-1},1\}$ and for any vertex $z\in V(P_{b}^{1})\setminus\{v_{1}\}$, it is the permutation of $\{p_{b},p_{2}\cdots,p_{k-1},1\}$. As $p_{b}\in[n]\setminus \{p_{1},p_{2}\cdots,p_{k-1},1\}$, then $p_{1}\neq p_{b}$ and $y\neq z$. Thus, the claim holds.

The proof of the Claim $1$ is complete.

\begin{clm}\label{clm2}
Let $X^{1}=\{u_{i}^{1}|u_{i}^{1}$ is the terminal vertex of the path $P_{i}^{1}$ for each $i\in\{2,3,\cdots,n\}\}$. Then the outside neighbours of vertices in $X^{1}$ belong to different copies of $B_{n-1,k-1}$, respectively.
\end{clm}

The proof of the Claim $2$. By Lemma~\ref{lem1}(2), the outside neighbours of vertices in $X^{1}$ are in $B_{n,k}^{2},B_{n,k}^{3},\cdots, B_{n,k}^{n}$, respectively. The proof of the Claim $2$ is complete.

\vskip0.3cm

Without loss of generality, suppose that the outside neighbour $(u_{i}^{1})^{\prime}$ of $u_{i}^{1}$ is in $B_{n,k}^{i}$ for each $i\in\{2,3,4,\cdots,n\}$. Otherwise, we can reorder the paths accordingly.

Similarly, let $v_{2}=p_{1}p_{2}p_{3}\cdots p_{k-1}2$, then there are $n-1$ paths $P_{1}^{2}, P_{3}^{2},\cdots, P_{n}^{2}$ starting at the vertex $v_{2}$ in $B_{n,k}^{2}$. Let $X^{2}=\{u_{1}^{2}, u_{3}^{2}, \cdots, u_{n}^{2}\}$ such that $u_{i}^{2}$ is the terminal vertex of the path $P_{i}^{2}$ and the outside neighbour $(u_{i}^{2})^{\prime}$ of $u_{i}^{2}$ is in $B_{n,k}^{i}$ for each $i\in\{1,3,4,\cdots,n\}$. In addition, there are $n-1$ paths $P_{1}^{3}, P_{2}^{3},\cdots, P_{n}^{3}$ starting at the vertex $v_{3}$ in $B_{n,k}^{3}$. Let $X^{3}=\{u_{1}^{3}, u_{2}^{3}, \cdots, u_{n}^{3}\}$ such that $u_{i}^{3}$ is the terminal vertex of the path $P_{i}^{3}$ and the outside neighbour $(u_{i}^{3})^{\prime}$ of $u_{i}^{3}$ is in $B_{n,k}^{i}$ for each $i\in\{1,2,4,\cdots,n\}$.

Obviously, the outside neighbour $(u_{1}^{3})^{\prime}$ of $u_{1}^{3}$ is in $B_{n,k}^{1}$ and the outside neighbour $(u_{2}^{3})^{\prime}$ of $u_{2}^{3}$ is in $B_{n,k}^{2}$. As $B_{n,k}^{1}$ is connected, there is a $((u_{1}^{3})^{\prime},v_{1})$-path $\widehat{P}_{1}$ in $B_{n,k}^{1}$. Let $t_{1}$ be the first vertex of the path $\widehat{P}_{1}$ which is in $\bigcup_{l\in\{2,3,\cdots,n\}}V(P_{l}^{1})$. Similarly, there is a $((u_{2}^{3})^{\prime},v_{2})$-path $\widehat{P}_{2}$ in $B_{n,k}^{2}$ as $B_{n,k}^{2}$ is connected. Let $t_{2}$ be the first vertex of the path $\widehat{P}_{2}$ which is in $\bigcup_{l\in\{1,3,\cdots,n\}}V(P_{l}^{2})$.
\begin{figure}[!ht]
\begin{center}
\vskip1cm
\includegraphics[scale=0.8]{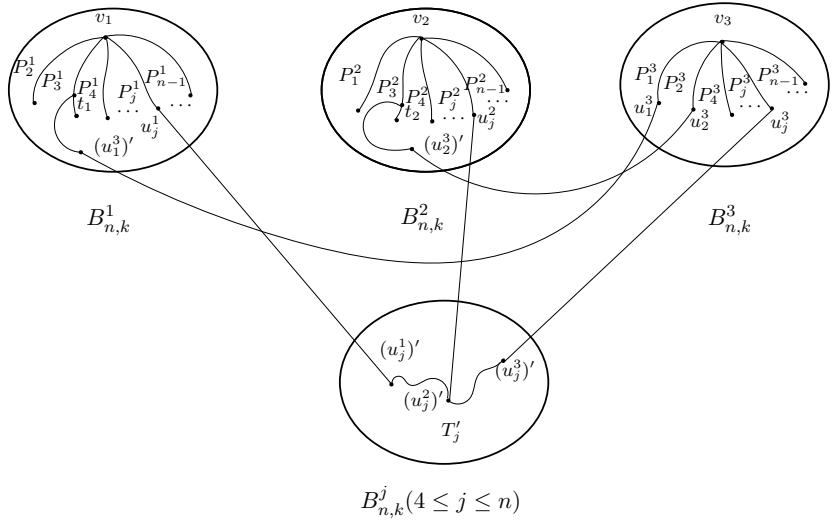}
\end{center}
\vskip0.5cm
\caption{The illustration of Subcase $3.1$ for $t_{1}\in V(P_{3}^{1})$ and $t_{2}\in V(P_{3}^{2})$ }\label{F2}
\end{figure}


To prove the result for $3\leq k\leq n-1$, the following two subcases are considered.

Subcase $3.1$. $t_{1}\in \bigcup_{l\in\{2,3\}} V(P_{l}^{1})$ and $t_{2}\in \bigcup_{l\in\{1,3\}} V(P_{l}^{2})$.

In this case, the induced subgraph $B_{n,k}[V(P_{1}^{3})\bigcup V(P_{2}^{1})\bigcup V(P_{3}^{1})\bigcup V(\widehat{P}_{1}[(u_{1}^{3})^{\prime}, t_{1}])]$ of $B_{n,k}$ contains a $(v_{3},v_{1})$-path, where $\widehat{P}_{1}[(u_{1}^{3})^{\prime}, t_{1}]$ is the subpath of $\widehat{P}_{1}$ starting at $(u_{1}^{3})^{\prime}$ and ending at $t_{1}$. Similarly, the induced subgraph $B_{n,k}[V(P_{2}^{3})\bigcup V(P_{1}^{2})\bigcup V(P_{3}^{2})\bigcup V(\widehat{P}_{2}[(u_{2}^{3})^{\prime}, t_{2}])]$ of $B_{n,k}$ contains a $(v_{3},v_{2})$-path, where $\widehat{P}_{2}[(u_{2}^{3})^{\prime}, t_{2}]$ is the subpath of $\widehat{P}_{2}$ starting at $(u_{2}^{3})^{\prime}$ and ending at $t_{2}$. The union of the $(v_{3}, v_{1})$-path and the $(v_{3}, v_{2})$-path forms a tree $T_{1}$ connecting $S$ in $B_{n,k}$. See Figure $2$.

In addition, as $(u_{j}^{1})^{\prime}, (u_{j}^{2})^{\prime}, (u_{j}^{3})^{\prime}\in V(B_{n,k}^{j})$ for each $j\in\{4,5,\cdots,n\}$ and $B_{n,k}^{j}$ is connected, there is a tree $T_{j}^{\prime}$ connecting $(u_{j}^{1})^{\prime}, (u_{j}^{2})^{\prime}$ and  $(u_{j}^{3})^{\prime}$ in $B_{n,k}^{j}$. Let $T_{j}=T_{j}^{\prime}\bigcup P_{j}^{1}\bigcup P_{j}^{2}\bigcup P_{j}^{3}\bigcup \\u_{j}^{1}(u_{j}^{1})^{\prime}\bigcup u_{j}^{2}(u_{j}^{2})^{\prime}\bigcup u_{j}^{3}(u_{j}^{3})^{\prime}$ for each $j\in\{4,5,\cdots,n\}$. Combining the trees $T_{j}s$ for $4\leq j\leq n$ and the tree $T_{1}$, and $n-2$ internally disjoint trees connecting $S$ in $B_{n,k}$ are obtained.

Subcase $3.2$. $t_{1}\in \bigcup_{l\in\{4,5,\cdots,n\}} V(P_{l}^{1})$ or $t_{2}\in \bigcup_{l\in\{4,5,\cdots,n\}} V(P_{l}^{2})$.

Without loss of generality, let $t_{1}\in V(P_{4}^{1})$.
Note that $v_{1}=p_{1}p_{2}\cdots p_{k-1}1$. By the assumption that the outside neighbor of the terminal vertex in $P_i^1$ is in
$B_{n,k}^i$ for $i\in \{2,3,\ldots, k\}$, one has that $v_{1}=23\cdots k1$. It implies that $p_i=i+1$ for $1\leq i\leq k-1$.

If $k\geq 4$, we obtain that $p_{k-1}\neq 2$ and $p_{3}=4$. For any vertex $v\in V(P_{4}^{1})$, $v$ is a permutation of $\{p_{1}, p_{2}, \cdots, p_{k-1},1\}$. Next, we consider the path $P_{2}^{1}$. Note that $u_{2}^{1}$ is the terminal vertex of $P_{2}^{1}$ and $u_{2}^{1}=p_{2}p_{3}\cdots p_{k-1}p_{1}1=34\cdots k21$. We can extend the path $P_{2}^{1}$ starting from $u_{2}^{1}$ as follows: $(3\underline{4}56\cdots k21)(35\underline{4}6\cdots k21)\cdots(35\cdots 26\underline k{4}1)$. Let $\widehat{u}_{2}^{1}=35\cdots 241$ and the extended path starting at $v_{1}$ and ending at $\widehat{u}_{2}^{1}$ be $\widehat{P}_{2}^{1}$. Then the outside neighbour of $\widehat{u}_{2}^{1}$ is in $B_{n,k}^{4}$.

If $k=3$ and $t_{1}\neq v_{1}$, then $v_1=231$ and $4\in[n]\setminus \{p_{1}, p_{2}, 1\}=\{4,5,\ldots, n\}$ and the vertex $t_{1}$ is a permutation of $\{4, p_{2},1\}=\{4,3,1\}$.
Note that $u_{2}^{1}=p_{2}21=321$. Now, we extend the path $P_{2}^{1}$ starting from $u_{2}^{1}$ to $\widehat{P}_{2}^{1}$, where $\widehat{P}_{2}^{1}={P}_{2}^{1}(421)(241)$. Let $\widehat{u}_{2}^{1}=241$. Now replacing $P_{2}^{1}$ with $\widehat{P}_{2}^{1}$,
The outside neighbor of terminal vertex $\widehat{u}_{2}^{1}$ of $\widehat{P}_{2}^{1}$ is in $B_{n,k}^{4}$.

Next, we prove the following claim.

\begin{clm}\label{clm3}
$V(\widehat{P}_{2}^{1})\bigcap V(P_{j}^{1})=\{v_{1}\}$ for each $j\in\{3,4,\cdots,n\}$ for $k\geq 3$.
\end{clm}

The proof of Claim $3$. For $k\geq 4$, we prove the result by contradiction. Suppose that there exists $l\in\{3,4,\cdots,n\}$ such that $|V(\widehat{P}_{2}^{1})\bigcap V(P_{l}^{1})|\geq 2$. Assume that $u\in V(\widehat{P}_{2}^{1})\bigcap V(P_{l}^{1})$ and $u\neq v_{1}$. Since $V(P_{2}^{1})\bigcap V(P_{l}^{1})=\{v_{1}\}, u\notin V(P_{2}^{1})$. Thus, $u\in V(\widehat{P}_{2}^{1})\setminus V(P_{2}^{1})$.

If $u\neq \widehat{u}_{2}^{1}$, then the element at position $k-1$ of $u$ is $2$. However, the element at position $k-1$ of each vertex in $V(P_{l}^{1})$ is $p_{k-1}$ or $k$. As $k\neq 2$ and $p_{k-1}\neq 2$, a contradiction.

Next, suppose $u=\widehat{u}_{2}^{1}$. The $k=4$ and $u=u_{4}^{1}$. However, the element at position $k-2$ of $u_{4}^{1}$ is $i_{k-1}$, a contradiction.

For $k=3$, let $x\in V(P_{m}^{1})$ for $4\leq m\leq n$, then it is a permutation of $\{m,3,1\}$. However, for any vertex $y\in V(\widehat{P}_{2}^{1}\setminus P_{2}^{1})$, it is a permutation of $\{4,2,1\}$. Thus, $x\neq y$.
The proof of the claim is complete.

\vskip0.3cm

Similarly, if $t_{2}\in V(P_{\ell }^{2})$ and $\ell \in\{4,5,\cdots,n\}$, we can extend the path $P_{2}^{2}$ to obtain the extended path, say $\widehat{P}_{2}^{2}$, such that the outside neighbour of the terminal vertex of the extended path $\widehat{P}_{2}^{2}$ is in $B_{n,k}^{\ell }$ and there is only one common vertex $v_{2}$ between the extended path and other paths $P_{j}s$ in $B_{n,k}^{2}$.

Since the induced subgraph $B_{n,k}[V(P_{1}^{3})\bigcup V(\widehat{P}_{1}[(u_{1}^{3})^{\prime}, t_{1}])\bigcup V(P_{4}^{1})]$ contains a $(v_{3}, v_{1})$-path, say $D_1$. Similarly, the induced subgraph $B_{n,k}[V(P_{2}^{3})\bigcup V(\widehat{P}_{2}[(u_{2}^{3})^{\prime}, t_{2}])\bigcup V(P_{4}^{1})]$ contains a $(v_{3}, v_{2})$-path, say $D_2$. A tree, say $T_{1}$, by combining $D_1$ and $D_2$ is obtained and  the tree $T_{1}$ connects $S$ in $B_{n,k}$.

 Similar as subcase $3.1$ just by replacing $P_{4}^{1}$ with $\widehat{P}_{2}^{1}$ as $t_{1}\in V(P_{4}^{1})$ or 
 replacing $P_{\ell}^{2}$ with $\widehat{P}_{2}^{2}$ if $t_{2}\in V(P_{\ell }^{2})$ for $\ell \in\{4,5,\cdots,n\}$,
 there is a tree $T_{j}$ connecting $S\bigcup V(B_{n,k}^{j})$ for each $j\in\{4,5,\cdots,n\}$ and $T_{j}s$ are internally disjoint $S$-trees. Combining the trees $T_{j}s$ for $4\leq j\leq n$ and the tree $T_{1}$, $n-2$ internally disjoint trees connecting $S$ in $B_{n,k}$ are obtained. Thus, the result is desired.

\hfill\qed

\section{Concluding remarks}
The generalized $k$-connectivity is a generalization of traditional connectivity. In this paper, we focus on the $(n,k)$-bubble-sort graph, denoted by $B_{n,k}$. We study the generalized $3$-connectivity of $B_{n,k}$ and show that $\kappa_{3}(B_{n,k})=n-2$ for $2\leq k\leq n-1$. So far, there are few results about the generalized $k$-connectivity for larger $k$. We are interested in this topic and we would like to study in this direction to show the corresponding results of $B_{n,k}$ for $k\geq 4$.

\section*{Acknowledgments}
This work was supported by the National Natural Science Foundation of China (No. 11731002), the Fundamental Research Funds for the Central Universities (No. 2016JBM071, 2016JBZ012) and the $111$ Project of China (B16002).

\end{document}